\pdfoutput=1
\documentclass[a4paper,10pt,oneside]{article}%
\usepackage[a4paper,textwidth=15.5cm,textheight=25.1cm]{geometry}
\usepackage{multicol}
\usepackage[utf8x]{inputenc}
\usepackage{amssymb}
\usepackage{amsmath}
\usepackage{synttree}
\usepackage{amsthm}
\usepackage{algorithm}
\usepackage[noend]{algorithmic}
\usepackage{todonotes}
\usepackage{amsfonts}
\usepackage{graphicx}
\usepackage{comment}
\usepackage{hyperref}
\usepackage{enumitem}
\usepackage{hyphenat}
\usepackage{tikz}%
\usepackage{authblk}
\usetikzlibrary{arrows}
\allowdisplaybreaks

\theoremstyle{definition}
\newtheorem{definition}{Definition}[section]
\theoremstyle{plain}
\newtheorem{theorem}[definition]{Theorem}
\newtheorem{proposition}[definition]{Proposition}

\theoremstyle{remark}
\newtheorem{remark}[definition]{Remark}
\newtheorem{example}[definition]{Example}

\title{\textsc{3D printing dimensional calibration shape: Clebsch Cubic}}

\date{}

\author[1]{Author Andr\'{e}. F. van der Merwe}
\author[2]{Janko B\"ohm}
\author[3]{Magdaleen S. Marais}
\affil[1]{Department of Industrial Engineering,
Stellenbosch University,
South Africa, andrevdm@sun.ac.za}
\affil[2]{Department of Mathematics, University of Kaiserslautern, Germany, boehm@mathematik.uni-kl.de }
\affil[3]{Department of Mathematics and Applied Mathematics, University of Pretoria, South Africa, magdaleen.marais@up.ac.za}

\begin{document}

\maketitle

\begin{abstract}
3D printing and other layer manufacturing processes are challenged by
dimensional accuracy. Several techniques are used to validate and calibrate
dimensional accuracy through the complete building envelope. The validation
process involves the growing and measuring of a shape with known parameters.
The measured result is compared with the intended digital model. Processes
with the risk of deformation after time or post processing may find this
technique beneficial. We propose to use objects from algebraic geometry as
test shapes. A cubic surface is given as the zero set of a $3$rd degree
polynomial with $3$ variables. A class of cubics in real $3$D space contains
exactly $27$ real lines. We provide a library for the computer algebra system
Singular which, from $6$ given points in the plane, constructs a cubic and the
lines on it. A surface shape derived from a cubic offers simplicity to the
dimensional comparison process, in that the straight lines and many other features can be analytically
determined and easily measured using non-digital equipment. For example,
the surface contains so-called Eckardt points, in each of which three of the lines intersect, and also other intersection points of pairs of lines. Distances
between these intersection points can easily be measured, since the points
are connected by straight lines. At all intersection
points of lines, angles can be verified. Hence, many features distributed
over the build volume are known analytically, and can be used for the validation process. Due to the thin shape geometry the
material required to produce an algebraic surface is minimal. This paper is
the first in a series that proposes the process chain to first define a cubic
with a configuration of lines in a given print volume and then to develop the
point cloud for the final manufacturing. Simple measuring techniques are recommended.

\end{abstract}

\section{Introduction}

This paper is the first in a series to investigate whether cubic surface shapes, and specifically the Clebsch cubic, can be used in 3D printing build volume accuracy. In this initial paper
the phases of development are proposed and the authors attempt to determine
the mathematical base for calculating with cubic surfaces. Various build
volumes, growing techniques and materials may require slight adjustments due to
its unique characteristics. However, the basic shape and mathematical approach
remains the same for all variants. The ultimate aim is to have a standard
Clebsch cubic shape which can be grown on any platform in any material, and in
any build volume size. The research phases proposed are an initial analytical
mathematical model, then an engine which converts from the analytical model to
a point cloud, then a digital domain simulated growth, followed by an actual
hardware printing phase, and lastly a reverse engineering phase. The initial
mathematical model is developed from ground rules to provide others the
fundamental information for parallel development. The input to the
mathematical model, based on the mathematical formulations found by Clebsch and others, is the
extent of the build volume. The open source computer algebra system
\textsc{Singular} is used for this conversion. The outputs from the
mathematical model is a three dimensional shape in analytical mathematical
formulation, the formulas for 27 straight lines, the coordinates of the points
where the lines cross, and the angles between the lines. After printing, the
line straightness is one indicator of the dimensional accuracy. Another
indicator is the angles between lines. The coordinate positioning of the cross
points, and the thickness of the cubic's vanes could also be measures. This
model is developed and suggested in the second part of this paper. The engine
which converts the analytical mathematical formulation to a printable point
cloud would typically be programmed in Matlab and later in C++. The inputs are
the three dimensional mathematical shape of the Clebsch cubic, and the
formulas of the 27 straight lines that we want to use as part of the
dimensional accuracy measurement. Note that the straight lines will
have to be highlighted in some way for the reverse engineering process to pick
it up. Several ways could be used to highlight the lines: generating
cylinders around the lines with diameter larger than the cubics vane
thickness, thinning the vanes along the lines, or perforating the
vanes along the lines, are examples. The output of this engine would be a
point cloud in an STL format or similar. In phase three of this
project we would attempt to compare the point cloud with the initial
analytical line formulas. This comparison can be done on a CAD platform, but
would typically be a manual process. Several alternatives of the previous
phases will be evaluated for accuracy. All work up to this point is in the
digital domain. This phase is to ensure accuracy and robustness of self
developed tools by comparison with trusted commercially available CAD
platforms. The output of this phase is a report which defines constraints and
extents within which these techniques are deemed accurate in the digital domain.
Phase four will see the growing of hard copies in various materials on various
platforms. This phase will report on any manufacturing issues on any of the
platforms using the extent of materials chosen. Phase five will reverse
engineer hardware shapes to compare with the initial intended analytical
shape. In this phase it will be determined to what extent the use of the $27$
straight lines and their angles is an indication of dimensional accuracy of
the process. This phase will seek to propose an economical method of measuring
dimensional accuracy of the complete building envelope. This paper starts with
the mathematical setup on which the description of the cubic and the lines are
based. Then a reference to the mathematical origins of the cubic surfaces is made, followed
by the derivation of the surface equation and the lines. Finally an example output from
the computer algebra library and explicit data for the Clebsch cubic is given.

\section{Algebraic varieties}

We first set the mathematical framework used to describe the cubic and
the lines on it. Let $K$ be either the real numbers $\mathbb{R}$ or the complex numbers $\mathbb{C}$. The set of
lines through the origin in $K^{n+1}$ is called \textbf{projective space} and
is denoted $\mathbb{P}_{K}^{n}$. We will write $(x_{0}:\ldots:x_{n})$ for the
line with direction $(x_{0},\ldots,x_{n})\neq0$. There is an inclusion of
usual $n$-space to projective space%
\[
K^{n}\longrightarrow\mathbb{P}_{K}^{n}\text{, }\hspace{3mm}(x_{1},\ldots
,x_{n})\longmapsto(1:x_{1}:\ldots:x_{n}).
\]
This map is referred to as an \textbf{affine chart}. The complement of the image is called the
plane at infinity (the horizon in a perspective drawing). An \textbf{algebraic
variety} $V(f_{1},\ldots,f_{r})\subset\mathbb{P}_{K}^{n}$ is the common zero
set of homogeneous polynomials $f_{i}\in K[x_{0},\ldots,x_{n}]$.

Algebraic varieties are studied in algebraic geometry, which forms a
central branch of classical mathematics. It has important applications, e.g.,
in cryptography, robotics, and computational biology.
 Algebraic varieties have the advantage over zero sets of 
non-polynomial equations that they can easily be handled by the means of computer algebra. For computing with polynomials we make use of the  open-source computer algebra system \textsc{Singular}
\cite{Singular}. Using projective
space in the development of the theory, avoids the problem that some
features of an algebraic variety (e.g. a line on it) may be contained in the plane at infinity.
For an introduction to algebraic geometry, computer algebra, and its
applications see, e.g. \cite{CLO}.

\section{Historic overview and derivation of the fundamental properties of
cubic hypersurfaces}

Starting in the second half of the $19$th century, Clebsch,\ Klein, Salmon,
Coble and many other mathematicians investigated cubic surfaces in
$\mathbb{P}_{\mathbb{C}}^{3}$, which are given by a single $3$rd degree polynomial. In $1849$, Arthur Cayley \cite{Cayley1849} and George Salmon
\cite{Salmon1849} found:

\begin{theorem}
Every smooth cubic surface in $\mathbb{P}_{\mathbb{C}}^{3}$ contains exactly
$27$ lines.
\end{theorem}

Here smooth means, that $C$ has in every point a well-defined tangent plane.
In algebraic geometry there is a process, called blowup, which replaces in a
variety a given point by a line and is a $1:1$ map everywhere else. In $1871$
Alfred Clebsch \cite{Clebsch} proved (see also \cite{Clebsch1866}):

\begin{theorem}
Every smooth cubic surface in $\mathbb{P}_{\mathbb{C}}^{3}$ is the blowup of
$\mathbb{P}_{\mathbb{C}}^{2}$ in $6$ points.
\end{theorem}

In the following, let $P_{1},\ldots,P_{6}\in\mathbb{P}_{K}^{2}$ be points in
general position, that is, no three are on a line and not all of them on a conic.

\begin{remark}
The homogeneous linear polynomial%
\[
l_{i,j}(t):=\det\left(  P_{i},P_{j},t\right)  :=\det\left(
\begin{array}
[c]{rrr}%
P_{i,0} & P_{j,0} & t_{0}\\
P_{i,1} & P_{j,1} & t_{1}\\
P_{i,2} & P_{j,2} & t_{2}%
\end{array}
\right)  \in K[t_{0},t_{1},t_{2}]
\]
defines in $\mathbb{P}_{K}^{2}$ the line through $P_{i}$ and $P_{j}$.
\end{remark}

\begin{proposition}
\cite{Coble} The blowup $C=C_{(P_{1},\ldots,P_{6})}$ of $\mathbb{P}_{K}^{2}$
in the points $P_{i}$ is the smallest algebraic variety (with respect to inclusion) containing the image
of%
\[%
\begin{tabular}
[c]{cccc}%
$\varphi_{(P_{1},\ldots,P_{6})}:$ & $\mathbb{P}_{K}^{2}\backslash
\{P_{1},\ldots P_{6}\}$ & $\longrightarrow$ & $\mathbb{P}_{K}^{5}$\\
& $(t_{0}:t_{1}:t_{2})$ & $\longmapsto$ & $\left(  \varphi_{0}(t):\ldots
:\varphi_{5}(t)\right)  $%
\end{tabular}
\ \
\]
(defined on $\mathbb{P}_{K}^{2}$ except at the points $P_{1},\ldots P_{6}$), where%
\[%
\begin{tabular}
[c]{l}%
$\varphi_{0}=l_{2,5}l_{1,3}l_{4,6}+l_{5,1}l_{4,2}l_{3,6}+l_{1,4}l_{3,5}%
l_{2,6}+l_{4,3}l_{2,1}l_{5,6}+l_{3,2}l_{5,4}l_{1,6}$\\
$\varphi_{1}=l_{5,3}l_{1,2}l_{4,6}+l_{1,4}l_{2,3}l_{5,6}+l_{2,5}l_{3,4}%
l_{1,6}+l_{3,1}l_{4,5}l_{2,6}+l_{4,2}l_{5,1}l_{3,6}$\\
$\varphi_{2}=l_{5,3}l_{4,1}l_{2,6}+l_{3,4}l_{2,5}l_{1,6}+l_{4,2}l_{1,3}%
l_{5,6}+l_{2,1}l_{5,4}l_{3,6}+l_{1,5}l_{3,2}l_{4,6}$\\
$\varphi_{3}=l_{4,5}l_{3,1}l_{2,6}+l_{5,3}l_{2,4}l_{1,6}+l_{4,1}l_{2,5}%
l_{3,6}+l_{3,2}l_{1,5}l_{4,6}+l_{2,1}l_{4,3}l_{5,6}$\\
$\varphi_{4}=l_{3,1}l_{2,4}l_{5,6}+l_{1,2}l_{5,3}l_{4,6}+l_{2,5}l_{4,1}%
l_{3,6}+l_{5,4}l_{3,2}l_{1,6}+l_{4,3}l_{1,5}l_{2,6}$\\
$\varphi_{5}=l_{4,2}l_{3,5}l_{1,6}+l_{2,3}l_{1,4}l_{5,6}+l_{3,1}l_{5,2}%
l_{4,6}+l_{1,5}l_{4,3}l_{2,6}+l_{5,4}l_{2,1}l_{3,6}$.
\end{tabular}
\ \
\]

\end{proposition}

\begin{remark}
\label{rmk clebsch}The \textbf{Clebsch cubic}, given in \cite[Ch. 16]%
{Clebsch}, is obtained by applying this construction to the points in general
position%
\[%
\begin{tabular}
[c]{lll}%
$P_{1}=(0:1:-g)$ & $P_{3}=(1:g:0)$ & $P_{5}=(0:1:g)$\\
$P_{2}=(g:0:1)$ & $P_{4}=(1:-g:0)$ & $P_{6}=(-g:0:1)$,
\end{tabular}
\]
where $g=\frac{1+\sqrt{5}}{2}$ is the golden ratio. These points correspond to
the diagonals in an icosahedron. The Clebsch cubic with $K=\mathbb{R}$
contains $27$ real lines.
\end{remark}

\begin{remark}\label{rmk triple det}
The number%
\[
\left\vert i,j;k,l;m,n\right\vert =\det\left(
\begin{array}
[c]{rr}%
\det\left(  P_{i},P_{j},P_{m}\right)  & \det\left(  P_{i},P_{j},P_{n}\right)
\\
\det\left(  P_{k},P_{l},P_{m}\right)  & \det\left(  P_{k},P_{l},P_{n}\right)
\end{array}
\right)  \text{,}%
\]
vanishes if the lines defined by $l_{i,j}(t)$, $l_{k,l}(t)$ and $l_{m,n}(t)$
in $\mathbb{P}_{K}^{2}$ meet in one point.
\end{remark}

\begin{theorem}
\cite{Coble,Cremona1878}\label{thm 37} Consider the skew-symmetric matrix%
\[
(A_{i,j})=\left(  {\footnotesize
\begin{array}
[c]{cccccc}%
0 & \left\vert 1,5;2,4;3,6\right\vert  & \left\vert 1,4;3,5;2,6\right\vert  &
\left\vert 1,2;4,3;5,6\right\vert  & \left\vert 2,3;4,5;1,6\right\vert  &
\left\vert 1,3;5,2;4,6\right\vert \\
& 0 & \left\vert 2,5;3,4;1,6\right\vert  & \left\vert 1,3;5,4;2,6\right\vert
& \left\vert 1,2;3,5;4,6\right\vert  & \left\vert 1,4;2,3;5,6\right\vert \\
&  & 0 & \left\vert 1,5;3,2;4,6\right\vert  & \left\vert
1,3;2,4;5,6\right\vert  & \left\vert 1,2;4,5;3,6\right\vert \\
&  &  & 0 & \left\vert 1,4;5,2;3,6\right\vert  & \left\vert
2,4;3,5;1,6\right\vert \\
& \overline{\phantom{xx}} &  &  & 0 & \left\vert 1,5;3,4;2,6\right\vert \\
&  &  &  &  & 0
\end{array}
}\right)\in K^{6\times 6}
\]
where the entries are defined as in Remark \ref{rmk triple det}, and write for the sum of the entries of the $i$-th row%
\[
a_{i}=\sum_{j=1}^{6}A_{i,j}\text{.}%
\]
Then $C$ is given by the equations
\begin{align*}
x_{0}^{3}+\ldots+x_{5}^{3}  &  =0\\
x_{0}+\ldots+x_{5}  &  =0\\
a_{0}\cdot x_{0}+\ldots+a_{5}\cdot x_{5}  &  =0\text{.}%
\end{align*}

\end{theorem}

\begin{remark}
Using the ordering of the $P_{i}$ from Remark \ref{rmk clebsch}, we obtain for
the Clebsch cubic surface $a_{0}=a_{1}=a_{2}=a_{3}=a_{4}=1$ and $a_{5}=-5$.
\end{remark}

\begin{remark}
For a subset $S\subset\mathbb{P}_{\mathbb{C}}^{n}$ we define $I(S)$ as the
ideal of all $f\in\mathbb{C}[x_{0},\ldots,x_{n}]$ with $f(x)=0$ for all $x\in
S$. So $V(I(S))$ is the smallest algebraic variety (with respect to inclusion) containing $S$. The ideal
generated by the $\varphi_{i}$ is
\[
\left\langle \varphi_{0},\ldots,\varphi_{5}\right\rangle =I(P_{1})\cap...\cap
I(P_{6})\text{.}%
\]
With the ring homomorphism%
\[%
\begin{tabular}
[c]{cccc}%
$\psi_{(P_{1},\ldots,P_{6})}:$ & $\mathbb{C}[x_{0},\ldots,x_{5}]$ &
$\longrightarrow$ & $\mathbb{C}[t_{0},t_{1},t_{2}]$\\
& $x_{i}$ & $\longmapsto$ & $\varphi_{i}$%
\end{tabular}
\ \ \
\]
we have
\[
I(C)=\ker\psi_{(P_{1},\ldots,P_{6})}=\left\langle x_{0}^{3}+...+x_{5}%
^{3},\text{ }x_{0}+...+x_{5},\text{ }a_{0}x_{0}+...+a_{5}x_{5}\right\rangle
\text{.}%
\]

\end{remark}

\begin{remark}
\label{rmk P3}Eliminating two variables by the two linear equations, $C$ can
be considered as a subset of $\mathbb{P}_{K}^{3}$.
\end{remark}

Note that a plane intersects $C$ in an irreducible plane cubic, a union of a
conic and a line, or in three lines.

\begin{definition}
A \textbf{tritangent plane }$H$\textbf{ to }$C$ is a plane, such that $H\cap
C$ consists out of three lines.
\end{definition}

\begin{remark}
A tritangent plane $H$ to $C$ is called \textbf{generic} if the three lines
pairwise intersect in three distinct points. Then $H$ is tangent to $C$ in
each of the three points.

If $H$ is not generic, then the three lines on $C$ intersect in a single
point. This point is called an \textbf{Eckardt point of }$C$.

Since in an Eckardt point the three lines are tangent to $C$, they are
coplanar, hence, lie on a tritangent plane. So, the Eckardt points are in
one-to-one correspondence to the non-generic tritangent planes.
\end{remark}

\begin{theorem}
\cite{Coble,Cayley1849} There are $45$ tritangent planes to $C$:

\begin{enumerate}
[leftmargin=0.8cm]

\item Of these, $15$ are given by the equations%
\[
x_{i}+x_{j}=0
\]
for $0\leq i<j\leq5$.

\item Write $M$ for the set of $2$-element subsets of $\{1,\ldots,6\}$, and
$S(M)$ for the set of permutations of $M$. The remaining $30$ tritangent
planes are then%
\[
\left(  m_{i,j}-d_{2}\right)  \cdot(x_{i}+x_{j})-(m_{k,l}+d_{2})\cdot
(x_{k}+x_{l})=0
\]
where
\[
\begin{tikzpicture}[|->,>=stealth',shorten >=1pt,auto,baseline={([yshift=-30pt]current bounding box.north)}]
\coordinate (a) at (0,.6);
\coordinate (b) at (-1,-.6);
\coordinate (c) at (1,-.6);
\node(1) at (a)  {\{i,j\}};
\node(2) at (b) {\{k,l\}};
\node(3) at (c)  {\{m,n\}};
\path[every node/.style={}]
(1) edge [bend right] node[left] {} (2)
(2) edge [bend right] node[left] {} (3)
(3) edge [bend right] node[right] {} (1);
\end{tikzpicture}\in S(M)
\]
is a $3$-cycle of pairwise disjoint elements of $M$,
\[
d_{2}={\scriptsize \det\left(
\begin{array}
[c]{rr}%
\det\left(  P_{3},P_{4},P_{1}\right)  \cdot\det\left(  P_{5},P_{6}%
,P_{1}\right)  & \det\left(  P_{5},P_{3},P_{1}\right)  \cdot\det\left(
P_{4},P_{6},P_{1}\right) \\
\det\left(  P_{3},P_{4},P_{2}\right)  \cdot\det\left(  P_{5},P_{6}%
,P_{2}\right)  & \det\left(  P_{5},P_{3},P_{2}\right)  \cdot\det\left(
P_{4},P_{6},P_{2}\right)
\end{array}
\right)  }%
\]
and%
\[
m_{i,j}=\sum_{s<t}a_{s}a_{t}+2(a_{i}^{2}+a_{j}^{2}+a_{i}a_{j})\text{,}%
\]
where $a_i$ is as defined in Theorem \ref{thm 37}.

\end{enumerate}
\end{theorem}

\begin{remark}
Possible numbers for Eckardt points are $1,2,3,4,6,9,10,18$. The Clebsch cubic
is the unique cubic with $10$ Eckardt points. The \textbf{Fermat cubic}
$V(x_{0}^{3}+\ldots+x_{3}^{3})$ is the unique cubic with the maximum possible
number of $18$ Eckardt points, however, only $3$ of the lines on the Fermat
cubic are defined over $\mathbb{R}$.
\end{remark}

\begin{remark}
Every line on $C$ lies on $5$ tritangent planes. Hence, any line on $C$ is the
intersection of the planes $x_{0}+\ldots+x_{5}=0$, $a_{0}x_{0}+...+a_{5}x_{5}=0$ and two
tritangent planes (see Remark \ref{rmk P3}).
\end{remark}

\begin{remark}
After permuting the coordinates we may assume that $a_{5}\neq0$. Then by
eliminating $x_{4}$ and $x_{5}$ via the two linear equations of $C$, we obtain
$C^{\prime}=V(F)\subset\mathbb{P}_{K}^{3}$ with a homogeneous cubic polynomial
$F\in K[x_{0},x_{1},x_{2},x_{3}]$.
\end{remark}

\begin{example}
The Clebsch Cubic is then given by
\[
F=x_{0}^{3}+x_{1}^{3}+x_{2}^{3}+x_{3}^{3}-\left(  x_{0}+x_{1}+x_{2}%
+x_{3}\right)  ^{3}\text{.}%
\]

\end{example}

\begin{remark}\label{rmk trafo}
For the Clebsch cubic, as well as cubics ``close'' to it in the sense of the
position of $P_{1},\ldots,P_{6}$, the transformation
\[%
\begin{tabular}
[c]{lll}%
$x_{0}=y_{0}-y_{3}-\sqrt{2}y_{1}$ &  & $x_{2}=y_{0}+y_{3}+\sqrt{2}y_{2}$\\
$x_{1}=y_{0}-y_{3}+\sqrt{2}y_{1}$ &  & $x_{3}=-y_{0}-y_{3}+\sqrt{2}y_{2}$%
\end{tabular}
\]
of the coordinate system with inverse%
\[%
\begin{tabular}
[c]{lll}%
$y_{0}=x_{0}+x_{1}+x_{2}-x_{3}$ &  & $y_{2}=\sqrt{2}(x_{2}+x_{3})$\\
$y_{1}=\sqrt{2}(-x_{0}+x_{1})$ &  & $y_{3}=-x_{0}-x_{1}+x_{2}-x_{3}$%
\end{tabular}
\
\]
achieves that all $27$ lines, for $K=\mathbb{R}$, are visible in the affine
chart%
\[
K^{3}\longrightarrow\mathbb{P}_{K}^{3}\text{, }\hspace{3mm}(y_{1},y_{2}%
,y_{3})\longmapsto(1:y_{1}:y_{2}:y_{3})\text{.}%
\]
Moreover, they all pass through a ball with radius $6$ around $0$. In the
affine chart we obtain a so-called affine cubic hypersurface $C^{\prime\prime
}\subset K^{3}$ given by a single, non-homogeneous $3$rd degree polynomial $f\in
K[y_{1},y_{2},y_{3}]$.
\end{remark}

\section{Implementation in Singular}

We have implemented the constructions above in the library \texttt{cubic.lib}
\cite{cubiclib} for the open-source computer algebra system \textsc{Singular}
\cite{Singular}. For an introduction to the language of \textsc{Singular} see
\cite{GP}. Specifically, from $6$ points in general position (with coordinates
in $\mathbb{Q}$ or an algebraic extension thereof), we give a function to
obtain the cubic $C\subset\mathbb{P}_{K}^{5}$, its projection $C^{\prime
}\subset\mathbb{P}_{K}^{3}$ and the affine cubic hypersurface $C^{\prime
\prime}\subset K^{3}$. Moreover, we compute the parametrizations
\[
\mathbb{P}_{K}^{2}\backslash\{P_{1},\ldots,P_{6}\}\longrightarrow
C\longrightarrow C^{\prime}%
\]
and an affine parametrization%
\[
\mathbb{P}_{K}^{2}\backslash V(\varphi_{0}+\varphi_{1}+\varphi_{2}-\varphi
_{3})\longrightarrow C^{\prime\prime}\text{.}%
\]
Finally, we compute the lines on $C,C^{\prime}$ and $C^{\prime\prime}$ in
implicit and parametric form, as well as the Eckardt points. We demonstrate
key parts of our library, considering the Clebsch cubic as an example:

\begin{example}
Our library can be loaded in \textsc{Singular} by:

\noindent\texttt{$>$ LIB "cubic.lib";}

\noindent We first create a polynomial ring in $4$ variables over the field
$\mathbb{Q}[\sqrt{5}]$:

\noindent\texttt{$>$ ring R = (0,a),(x0,x1,x2,x3),dp;}

\noindent\texttt{$>$ minpoly = a\symbol{94}2-5;}

\noindent We specify a list $P$ with the points $P_{1},\ldots,P_{6}$:

\noindent\texttt{$>$ number g = (1 + a)/2;}

\noindent\texttt{$>$ list P = vector(0,1,-g), vector(g,0,1), vector(1,g,0),
vector(1,-g,0),}

\noindent\hspace{2.35cm}\texttt{vector(0,1,g), vector(-g,0,1);}

\noindent We compute the equation of $C^{\prime}$:

\noindent\texttt{$>$ poly f =\ cubic(P);}

\noindent\texttt{$>$ f;}

\noindent\texttt{-3*x0\symbol{94}2*x1-3*x0*x1\symbol{94}2-3*x0\symbol{94}%
2*x2-6*x0*x1*x2-3*x1\symbol{94}2*x2-3*x0*x2\symbol{94}2-3*x1*x2\symbol{94}2}

\noindent\texttt{-3*x0\symbol{94}2*x3-6*x0*x1*x3-3*x1\symbol{94}%
2*x3-6*x0*x2*x3-6*x1*x2*x3-3*x2\symbol{94}2*x3-3*x0*x3\symbol{94}2}

\noindent\texttt{-3*x1*x3\symbol{94}2-3*x2*x3\symbol{94}2}

\noindent The following command returns a list of all lines on $C^{\prime}$,
each specified by $2$ linear equations:

\noindent\texttt{$>$ list L = lines(P);}

\noindent\texttt{$>$ L[1];}

\noindent\texttt{\_[1] = x0 + x1}

\noindent\texttt{\_[2] = x2 + x3}

\noindent We compute a list of Eckardt points, each specified by $3$ linear equations:

\noindent\texttt{$>$ list E = EckardtPoints(P);}

\noindent\texttt{$>$ E[1];}

\noindent\texttt{\_[1] = x0}

\noindent\texttt{\_[2] = x1}

\noindent\texttt{\_[3] = x2 + x3}

\noindent Note that, by the commands \texttt{affineCubic},
\texttt{affineLines} and \texttt{affineEckardtPoints}, one can also obtain the
affine cubic $C^{\prime\prime}$ and the corresponding lines and Eckardt
points, respectively. Moreover, the functions \texttt{paramLines} and
\texttt{affineParamLines} compute parametrizations of the lines on $C^{\prime
}$ and $C^{\prime\prime}$, respectively.

If, in addition to \textsc{Singular}, the program \textsc{Surf} \cite{surf} is
installed, $C^{\prime\prime}$ can be visualized by:

\noindent\texttt{$>$ LIB "surf.lib";}

\noindent\texttt{$>$ plot(affineCubic(P));}

\noindent It also can plot hyperplane sections of a surface. Hence, we can
visualize the lines on the cubic by intersecting with tritangent planes, see
Figure \ref{fig cubic}.\begin{figure}[h]
\begin{center}
\includegraphics[
natheight=15.000300in,
natwidth=15.000300in,
height=4.2117in,
width=4.2117in
]{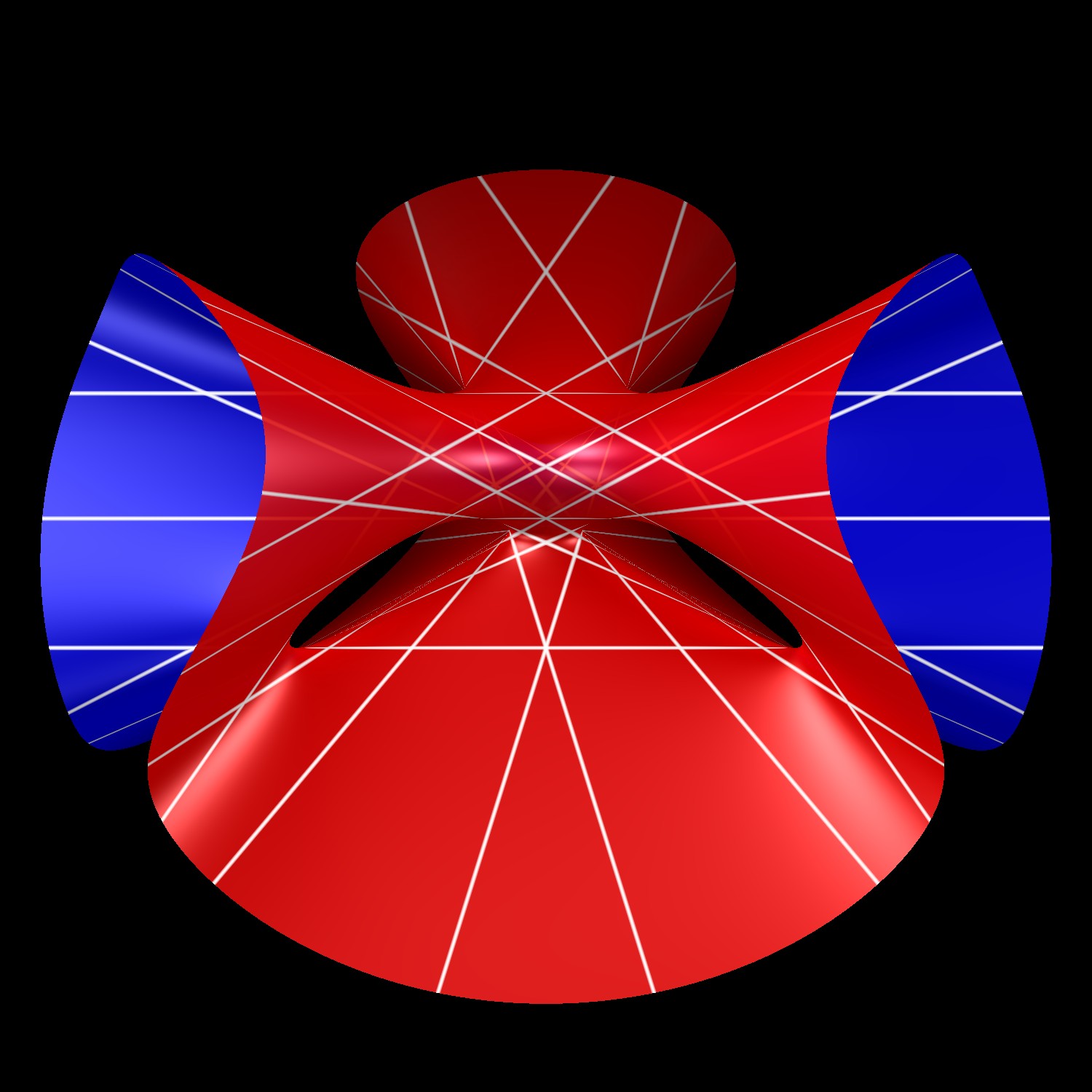}
\end{center}
\caption{Lines on the Clebsch cubic}%
\label{fig cubic}%
\end{figure}
\end{example}

\section{Explicit data for the Clebsch cubic}

Our program returns the equation for the cubic and derived information for any
6 points in general position in the projective plane (corresponding to lines
in 3D space through the origin).

In this section we give the explicit data required for the dimensional
comparison process for the Clebsch cubic (i.e. for the choice of the six
points given by the diagonals in an icosahedron). In the following let $a=\sqrt{5}$ and 
$c=\sqrt{2}$.
The cubic $C^{\prime\prime}$ is the zero set in $K^{3}$ of the
equation
\begin{align*}
& 2 c y_2^3+2 y_1^2 y_3-8 y_2^2 y_3+3 c y_2 y_3^2-y_3^3-2 y_1^2+8 y_2^2-10 c y_2 y_3+3 y_3^2+3 c y_2-3 y_3+1 = 0\text{.}
\end{align*}
The lines on $C^{\prime\prime}$ given in implicit form (by two linear equations each) as well as their parametrizations (specified as maps $K\rightarrow K^3, s \mapsto (\psi_1 (s),\psi_2 (s),\psi_3 (s))$) are specified in Table \ref{tab lines}.
\begin{table}[h]
\label{tab lines}
${\tiny
\def\arraystretch{1.5}
\begin{array}
[c]{c|c}%
\text{\begin{normalsize}Implicit\end{normalsize}} & \text{\begin{normalsize}Parametric\end{normalsize}}\\
\hline\\
\left\langle y_{3}-1,y_{2}\right\rangle  & \left(
s
,
0
,
1
\right)\\
\left\langle y_{3}-1,y_{2}-\frac{2}{c}\right\rangle  & \left(
s
,
\frac{2}{c}
,
1
\right) \\
\left\langle y_{3}-1,y_{2}+\frac{2}{c}\right\rangle  & \left(
s
,
-\frac{2}{c}
,
1
\right) \\
\left\langle y_{2}-\frac{1}{c}\cdot y_{3}+\frac{1}{c},y_{1}-\frac{1}{c}\cdot
y_{3}-\frac{1}{c}\right\rangle  & \left(
\frac{1}{c}\cdot s+\frac{1}{c}
,
\frac{1}{c}\cdot s-\frac{1}{c}
,
s
\right)\\
\left\langle y_{2}-\frac{1}{c}\cdot y_{3}-\frac{1}{c},y_{1}-\frac{1}{c}\cdot
y_{3}-\frac{3}{c}\right\rangle  & \left(
\frac{1}{c}\cdot s+\frac{3}{c}
,
\frac{1}{c}\cdot s+\frac{1}{c}
,
s
\right) \\
\left\langle y_{2}-\frac{1}{c}\cdot y_{3}+\frac{3}{c},y_{1}-\frac{1}{c}\cdot
y_{3}+\frac{1}{c}\right\rangle  & \left(
\frac{1}{c}\cdot s-\frac{1}{c}
,
\frac{1}{c}\cdot s-\frac{3}{c}
,
s
\right)\\
\left\langle y_{2}-\frac{1}{c}\cdot y_{3}+\frac{1}{c},y_{1}+\frac{1}{c}\cdot
y_{3}+\frac{1}{c}\right\rangle  & \left(
-\frac{1}{c}\cdot s-\frac{1}{c}
,
\frac{1}{c}\cdot s-\frac{1}{c}
,
s
\right) \\
\left\langle y_{2}-\frac{1}{c}\cdot y_{3}-\frac{1}{c},y_{1}+\frac{1}{c}\cdot
y_{3}+\frac{3}{c}\right\rangle  & \left(
-\frac{1}{c}\cdot s-\frac{3}{c}
,
\frac{1}{c}\cdot s+\frac{1}{c}
,
s
\right) \\
\left\langle y_{2}-\frac{1}{c}\cdot y_{3}+\frac{3}{c},y_{1}+\frac{1}{c}\cdot
y_{3}-\frac{1}{c}\right\rangle  & \left(
-\frac{1}{c}\cdot s+\frac{1}{c}
,
\frac{1}{c}\cdot s-\frac{3}{c}
,
s
\right) \\
\left\langle y_{2}+\frac{1}{c}\cdot y_{3}+\frac{1}{c},y_{1}+\frac{3}{c}\cdot
y_{3}+\frac{1}{c}\right\rangle  & \left(
-\frac{3}{c}\cdot s-\frac{1}{c}
,
-\frac{1}{c}\cdot s-\frac{1}{c}
,
s
\right) \\
\left\langle y_{2}-\frac{3}{c}\cdot y_{3}+\frac{1}{c},y_{1}-\frac{1}{c}\cdot
y_{3}+\frac{1}{c}\right\rangle  & \left(
\frac{1}{c}\cdot s-\frac{1}{c}
,
\frac{3}{c}\cdot s-\frac{1}{c}
,
s
\right) \\
\left\langle y_{2}+\frac{1}{c}\cdot y_{3}+\frac{1}{c},y_{1}-\frac{3}{c}\cdot
y_{3}-\frac{1}{c}\right\rangle  & \left(
\frac{3}{c}\cdot s+\frac{1}{c}
,
-\frac{1}{c}\cdot s-\frac{1}{c}
,
s
\right) \\
\left\langle y_{2}-\frac{3}{c}\cdot y_{3}+\frac{1}{c},y_{1}+\frac{1}{c}\cdot
y_{3}-\frac{1}{c}\right\rangle  & \left(
-\frac{1}{c}\cdot s+\frac{1}{c}
,
\frac{3}{c}\cdot s-\frac{1}{c}
,
s
\right) \\
\left\langle y_{2},y_{1}-\frac{1}{c}\cdot y_{3}+\frac{1}{c}\right\rangle  & \left(
\frac{1}{c}\cdot s-\frac{1}{c}
,
0
,
s
\right) \\
\left\langle y_{2},y_{1}+\frac{1}{c}\cdot y_{3}-\frac{1}{c}\right\rangle  & \left(
-\frac{1}{c}\cdot s+\frac{1}{c}
,
0
,
s
\right)\\
\left\langle y_{2}-\frac{1}{ac+2c}\cdot y_{3}+\frac{5}{ac},y_{1}+(\frac{3}%
{2}ac-3c)\cdot y_{3}-(\frac{1}{2}ac-2c)\right\rangle  & \left(
-(\frac{3}{2}ac-3c)\cdot s+(\frac{1}{2}ac-2c)
,
\frac{1}{ac+2c}\cdot s-\frac{1}{2}ac
,
s
\right) \\
\left\langle y_{2}-\frac{5}{ac}\cdot y_{3}+\frac{1}{ac+2c},y_{1}-\frac
{11}{ac+4c}\cdot y_{3}-\frac{3}{ac+2c}\right\rangle  & \left(
\frac{11}{ac+4c}\cdot s+\frac{3}{ac+2c}
,
\frac{5}{ac}\cdot s-\frac{1}{ac+2c}
,
s
\right) \\
\left\langle y_{2}+\frac{1}{ac-2c}\cdot y_{3}-\frac{5}{ac},y_{1}-(\frac{3}%
{2}ac+3c)\cdot y_{3}+(\frac{1}{2}ac+2c)\right\rangle  & \left(
(\frac{3}{2}ac+3c)\cdot s-(\frac{1}{2}ac+2c)
,
-\frac{1}{ac-2c}\cdot s+\frac{1}{2}ac
,
s
\right) \\
\left\langle y_{2}+\frac{5}{ac}\cdot y_{3}-\frac{1}{ac-2c},y_{1}+\frac
{11}{ac-4c}\cdot y_{3}+\frac{3}{ac-2c}\right\rangle  & \left(
-\frac{11}{ac-4c}\cdot s-\frac{3}{ac-2c}
,
-\frac{5}{ac}\cdot s+\frac{1}{ac-2c}
,
s
\right)\\
\left\langle y_{2}+\frac{2}{ac-3c}\cdot y_{3}-\frac{2}{ac-3c},y_{1}+(\frac
{1}{4}ac+\frac{1}{4}c)\cdot y_{3}+(\frac{1}{4}ac+\frac{1}{4}c)\right\rangle  &
\left(
-(\frac{1}{4}ac+\frac{1}{4}c)\cdot s-(\frac{1}{4}ac+\frac{1}{4}c)
,
-\frac{2}{ac-3c}\cdot s-(\frac{1}{4}ac+\frac{3}{4}c)
,
s
\right)\\
\left\langle y_{2}+\frac{1}{ac-2c}\cdot y_{3}-\frac{5}{ac},y_{1}+(\frac{3}%
{2}ac+3c)\cdot y_{3}-(\frac{1}{2}ac+2c)\right\rangle  & \left(
-(\frac{3}{2}ac+3c)\cdot s+(\frac{1}{2}ac+2c)
,
-\frac{1}{ac-2c}\cdot s+\frac{1}{2}ac
,
s
\right) \\
\left\langle y_{2}-\frac{2}{ac+3c}\cdot y_{3}+\frac{2}{ac+3c},y_{1}-(\frac
{1}{4}ac-\frac{1}{4}c)\cdot y_{3}-(\frac{1}{4}ac-\frac{1}{4}c)\right\rangle  &
\left(
(\frac{1}{4}ac-\frac{1}{4}c)\cdot s+(\frac{1}{4}ac-\frac{1}{4}c)
,
\frac{2}{ac+3c}\cdot s+(\frac{1}{4}ac-\frac{3}{4}c)
,
s
\right)
\\
\left\langle y_{2}-\frac{1}{ac+2c}\cdot y_{3}+\frac{5}{ac},y_{1}-(\frac{3}%
{2}ac-3c)\cdot y_{3}+(\frac{1}{2}ac-2c)\right\rangle  & \left(
(\frac{3}{2}ac-3c)\cdot s-(\frac{1}{2}ac-2c)
,
\frac{1}{ac+2c}\cdot s-\frac{1}{2}ac
,
s
\right)\\
\left\langle y_{2}+\frac{2}{ac-3c}\cdot y_{3}-\frac{2}{ac-3c},y_{1}-(\frac
{1}{4}ac+\frac{1}{4}c)\cdot y_{3}-(\frac{1}{4}ac+\frac{1}{4}c)\right\rangle  &
\left(
(\frac{1}{4}ac+\frac{1}{4}c)\cdot s+(\frac{1}{4}ac+\frac{1}{4}c)
,
-\frac{2}{ac-3c}\cdot s-(\frac{1}{4}ac+\frac{3}{4}c)
,
s
\right)\\
\left\langle y_{2}+\frac{5}{ac}\cdot y_{3}-\frac{1}{ac-2c},y_{1}-\frac
{11}{ac-4c}\cdot y_{3}-\frac{3}{ac-2c}\right\rangle  & \left(
\frac{11}{ac-4c}\cdot s+\frac{3}{ac-2c}
,
-\frac{5}{ac}\cdot s+\frac{1}{ac-2c}
,
s
\right)\\
\left\langle y_{2}-\frac{2}{ac+3c}\cdot y_{3}+\frac{2}{ac+3c},y_{1}+(\frac
{1}{4}ac-\frac{1}{4}c)\cdot y_{3}+(\frac{1}{4}ac-\frac{1}{4}c)\right\rangle  & \left(
-(\frac{1}{4}ac-\frac{1}{4}c)\cdot s-(\frac{1}{4}ac-\frac{1}{4}c)
,
\frac{2}{ac+3c}\cdot s+(\frac{1}{4}ac-\frac{3}{4}c)
,
s
\right)
\\
\left\langle y_{2}-\frac{5}{ac}\cdot y_{3}+\frac{1}{ac+2c},y_{1}+\frac
{11}{ac+4c}\cdot y_{3}+\frac{3}{ac+2c}\right\rangle  & \left(
-\frac{11}{ac+4c}\cdot s-\frac{3}{ac+2c}
,
\frac{5}{ac}\cdot s-\frac{1}{ac+2c}
,
s
\right) \\
\\
\hline
\end{array}
}$
\medskip
\caption{Lines on the Clebsch cubic in implicit and parametric form}
\end{table}
The $10$ Eckardt points on $C^{\prime}$ have projective coordinates%
\[%
\begin{tabular}
[c]{ccccc}%
$(-1:1:0:0)$ & $(-1:0:1:0)$ & $(0:-1:1:0)$ & $(-1:0:0:1)$ & $(0:-1:0:1)$\\
$(0:0:-1:1)$ & $(1:0:0:0)$ & $(0:1:0:0)$ & $(0:0:1:0)$ & $(0:0:0:1)$.%
\end{tabular}
\
\]
Hence (after applying the transformation of Remark \ref{rmk trafo} and passing to the affine
chart), the cubic $C^{\prime\prime}$ contains $7$ of them with affine
coordinates%
\[%
\begin{tabular}
[c]{ccccccc}%
$(-\frac{1}{c},-\frac{1}{c},0)$ & $(\frac{1}{c},-\frac{1}{c},0)$ & $(0,0,1)$ &
$(-\frac{2}{c},0,-1)$ & $(\frac{2}{c},0,-1)$ &
$(0,\frac{2}{c},1)$ & $(0,-\frac{2}{c},1)$, 
\end{tabular}
\
\]
the remaining three of them lying at the plane at infinity $y_{0}=0$ with projective
coordinates $(0:1:0:0)$, $(0:1:1:c)$, $(0:-1:1:c)$.

Using the implicit equations of the lines, it is easy to determine, which lines pass through which Eckardt points. The angle $\alpha$ beween two such lines with direction vectors $v_1$ and $v_2$ can be calculated
via the well-known formula $\cos(\alpha) = \frac{\left\langle v_1, v_2 \right\rangle}{ \left\lVert v_1 \right\rVert \cdot \left\lVert v_2 \right \rVert }$.

\begin{remark}
For the given data a suitable print volume is the cube
\[
[-6,6]\times[-6,6]\times[-6,6]\text{.}
\]
However, for practical purposes the build volume may be a rectangular cuboid
\[
[-r_{1},r_{1}]\times[-r_{2},r_{2}]\times[-r_{3},r_{3}]\text{.}
\]
From the above data, a suitable test shape is then obtaind by applying the
substitution
\[%
\begin{tabular}
[c]{lllll}%
$y_{1}\mapsto x\cdot\frac{6}{r_{1}}$ &  & $y_{2}\mapsto y_{2}\cdot\frac
{6}{r_{2}}$ &  & $y_{3}\mapsto y_{3}\cdot\frac{6}{r_{3}}$%
\end{tabular}
\]
to the implicit equations. Correspondingly, a point $(y_{1},y_{2},y_{3}%
)\in\mathbb{R}^{3}$ is transformed to the point $(y_{1}\cdot\frac{r_{1}}%
{6},\,y_{2}\cdot\frac{r_{2}}{6},\,y_{3}\cdot\frac{r_{3}}{6})$.
\end{remark}

\section{Conclusion}
A general cubic surface contains $27$ lines. For the Clebsch cubic all lines
are on its shaped surface defined over the real numbers. These straight lines are
determined analytically and are proposed for use to measure the build accuracy of a 3D printing process using even non digital techniques. The straight lines extend to the build volume
perimeter and are structurally supported by the cubic surface. It is proposed that the lines
are ``highlighted'' to stand out from the cubic for easy identification and physical measurement. The development procedure is outlined in this paper and is proposed as phases of the
project. The first phase, which is descibed by this paper, is the mathematical description of cubic surfaces and the formulas of its straight lines. Future papers will present the results
from the next phases that would firstly digitally test the accuracy and lastly produce physical parts for reverse engineering quality control purposes.
Moreover, we will describe a method to specify $6$ lines and then derive
a cubic test shape containing these given lines.

\bibliographystyle{plain}
\bibliography{bibcubic}

\end{document}